\documentclass[letterpaper,9pt]{IEEEtran}
\usepackage{cite}
\usepackage{amsmath,amssymb,amsfonts}
\usepackage{algorithmic}
\usepackage{graphicx}
\usepackage{textcomp}
\usepackage{xcolor}
\usepackage{mathrsfs}

\usepackage{amsfonts}
\usepackage{cite}

\usepackage{algorithm}
\usepackage{algorithmic}

\usepackage{indentfirst}

\usepackage{flushend}

\usepackage{multirow}
\usepackage{float}
\usepackage{algorithm}
\usepackage{algorithmic}
\usepackage{nomencl}

\usepackage{etoolbox}

\allowdisplaybreaks

\begin{document}

\begin{@twocolumnfalse}

    \title{
            \begin{flushleft}
            \Large \textbf{Preprint}
        \end{flushleft}
    \vspace{0.5cm}
        \large \textbf{Universal Stabilization with Six Records}}
    \author{Senhan Yao
\thanks{Corresponding author: Senhan Yao.}
\thanks{Senhan Yao is with the Key Laboratory of Systems and
Control, AMSS, Chinese Academy of Sciences, Beijing 100190,
China, and also with the School of Mathematical Sciences, University
of Chinese Academy of Sciences, Beijing 100049, China (e-mail: yaosenhan@amss.ac.cn).}
}
    \maketitle
\end{@twocolumnfalse}

\vspace{-1cm}\noindent \textbf{Abstract.} This letter shows that a single feedback controller needs only six retained observation records to handle the largest universally stabilizable uncertainty range for the scalar discrete-time benchmark model of Xie and Guo, when uncertainty is measured by the radius of a generalized Lipschitz ball. Existing constructions can require arbitrarily many records along bounded trajectories. Our controller retains at most six genuine state–response pairs and five auxiliary real scalars. Each update requires a uniformly bounded number of elementary operations. No knowledge of the plant, its generalized Lipschitz constants, or the disturbance bound is required.

Control engineers want to keep plants stable. Yet plant dynamics are rarely known in full. The available knowledge may constrain the possible dynamics to a model class. This class may include open-loop unstable plants. To ensure stability, we seek a controller that stabilizes the entire class. We call it a universal stabilizer over this class \cite{ByrnesWillems1984}. Such controllers are also called robustly stabilizing controllers or simultaneous stabilizers \cite{GloverMcFarlane1989, VidyasagarViswanadham1982}. The design of such controllers has been extensively studied. It remains a central goal of control theory.

Much of this work follows a common approach. Control theorists first construct a controller. They then prove that it universally stabilizes an uncertainty class. Identifying the largest class stabilized by the controller is usually difficult. The class is chosen to make the analysis tractable. Thus, the resulting guarantee is therefore often conservative: the controller may stabilize a substantially larger class than the analysis establishes.

Such conservatism leaves a fundamental question open. How far can an uncertainty class be enlarged while still admitting a universal stabilizer? Xie and Guo take a pioneering step in this direction \cite{XieGuo2000}. They seek a benchmark model that isolates uncertainty from limitations of the control structure. This allows uncertainty to be measured directly through the size of a class of possible plant functions. To this end, they study the scalar discrete-time system
\begin{equation}\label{eq:plant}
 x_{t+1}=f(x_t)+u_t+w_{t+1},
 \qquad t=0,1,\ldots,
\end{equation}
where \(x_t,u_t\in\mathbb{R}\), the plant \(f: \mathbb{R} \to \mathbb{R}\) is fixed and unknown, and \(\{w_t\}_{t \geq 1}\) is an arbitrary bounded disturbance sequence with an unknown bound \(W\). They allow all causal feedback laws \(u_t=h_t(x_0,\ldots,x_t)\), where each \(h_t:\mathbb R^{t+1}\to\mathbb R\) may be nonlinear and time-varying. A universal stabilizer must keep the state and control bounded for every plant in the class, every initial state, and every bounded disturbance sequence, i.e.
\begin{equation}
\forall f\in\mathcal C,\quad
\forall x_0\in\mathbb R,\quad
\forall w\in\ell_\infty(\mathbb N;\mathbb R),\quad
\sup_{t\ge0}\bigl(|x_t|+|u_t|\bigr)<\infty.
\end{equation}

To measure the size of an uncertainty class, they introduce the generalized Lipschitz norm
\begin{equation}\label{eq:GL}
\|f\|_{\mathrm{GL}}
:=
\lim_{\alpha\to\infty}
\sup_{x,y\in\mathbb R}
\frac{|f(x)-f(y)|}{|x-y|+\alpha}.
\end{equation}
They describe uncertainty through balls in this norm. For a ball centered at a known, locally bounded nominal plant, the nominal dynamics can be canceled through the additive control input. This transformation preserves causality and boundedness of both the state and the control. It therefore suffices to consider balls centered at zero. Among these balls, Xie and Guo identify the largest one that admits a universal stabilizer, which we call the Xie--Guo positive regime. This largest ball is open and has radius \(L^\star = 3/2 + \sqrt{2}\) known as the Xie--Guo constant \cite{Sokolov2016Automatica}. They construct a universal stabilizer for this ball, which we call the Xie--Guo stabilizer. They also prove that no larger ball admits a universal stabilizer. Thus, the Xie--Guo positive regime reflects the maximum uncertainty that can be dealt with by feedback mechanism.

Further developments in this field either enlarge the uncertainty class by admitting additional plant dynamics or consider systems with more complicated feedback structures \cite{Sokolov2016Automatica,ZhangGuo2002,XueGuo2002,Sokolov2003a,Ma2008a,HuangGuo2012,RenChengGuo2014,Sokolov2017,Sokolov2018ARC}. Each extension makes controller design and stability analysis more complicated. Within the simple feedback structure considered here, it remains unclear whether successive enlargements can reach a maximal class, i.e. a universally stabilizable class that cannot be enlarged further. This is more than a gap in the literature, because proving that such a class exists is itself mathematically difficult. If a maximal class does exist, designing a computationally efficient universal stabilizer for it presents another challenge. We can hardly expect to address this problem well without first developing an efficient universal stabilizer for the benchmark class already at hand.

We start with the Xie--Guo stabilizer \cite{XieGuo2000}. We first observe that each state transition provides a state-response pair
\begin{equation}\label{eq:visit_pair}
(x_i,z_i),\qquad z_i:=x_{i+1}-u_i=f(x_i)+w_{i+1}.
\end{equation}
Here we call the noisy observation of the plant dynamics a response. The Xie--Guo stabilizer stores these pairs. At each step \(t\), it finds the stored state \(x_{i_t}\) closest to the current state \(x_t\). It uses the associated response \(z_{i_t}\) to estimate the value of the plant dynamics at the current state \(f(x_t)\). The control input essentially subtracts this estimate from the midpoint of the observed state range, which is known as the midpoint stabilization law:
\begin{equation}\label{eq:control}
u_t=-z_{i_t}
+\frac12\left(
m_t+M_t
\right),
\end{equation}
where
\begin{equation}\label{eq:state_extrema}
m_t:=\min_{0\le i\le t}x_i,
\qquad
M_t:=\max_{0\le i\le t}x_i.
\end{equation}
are the smallest and largest states observed up to time \(t\). Its direct implementation stores and searches all past pairs. The number of stored pairs grows with time. Any fixed memory budget will therefore eventually be exhausted. A direct search also examines an increasing number of records at each step. With a fixed sampling period, the controller must finish this search before the next input is due. The implementation provides no bound on the work required per step, so this deadline cannot be guaranteed throughout continued operation.

Following a suggestion of Xie and Guo \cite[Remark~2.5]{XieGuo2000}, Sokolov reduces storage by selecting which observations to retain \cite[Sec.~III]{Sokolov2003a}. He fixes a spatial resolution and divides the state space into intervals of that width. He retains one observation record per visited interval. It is obvious that a bounded trajectory visits only finitely many intervals. Storage therefore remains finite along each trajectory. Sokolov discusses how resolution and proximity to the critical Lipschitz threshold \(L^\star\) affect record requirements \cite[Remarks~2--3, p.~609]{Sokolov2003a}. However, he does not quantify the claimed growth in the required record count. In fact, for fixed resolution, we can prove a record-count lower bound of order \((L^\star-L)^{-1/2}\), where \(L<L^\star\) denotes the generalized Lipschitz constant of the plant \(f\). In particular, record count can be arbitrarily large on the Xie--Guo positive regime.

\textbf{Record-Count Lower Bound for Sokolov's Scheme.} Fix \(0<\varepsilon<1/2\) and take \(f_L(x)=Lx+1\), with \(2\le L<L^\star\), \(x_0=u_0=0\), and no disturbance. Use Sokolov's Scheme and Xie--Guo stabilizer \eqref{eq:control} with zero nominal plant. When the current state is farther than \(\varepsilon\) from all stored states, this controller subtracts the estimated dynamics from the average of the smallest and largest observed states. Along an increasing positive trajectory, these extremes are \(0\) and \(x_t\), and the nearest stored state is \(x_{t-1}\). Thus \(u_t=-f_L(x_{t-1})+x_t/2\), giving \(x_{t+1}=(L+1/2)x_t-Lx_{t-1}\). Set \(\theta_L=\arccos((L+1/2)/(2\sqrt L))\) and \(n_L=\lfloor\pi/(3\theta_L)\rfloor\). The recurrence gives \(x_t=L^{(t-1)/2}\sin(t\theta_L)/\sin\theta_L\), whose increments satisfy \(x_t-x_{t-1}\ge1/2>\varepsilon\) for \(1\le t\le n_L+1\). This verifies the construction inductively and forces at least \(n_L\) completed records. Since \(\theta_L\sim\sqrt{(L^\star-L)/(\sqrt2L^\star)}\), the record count is \(\Omega((L^\star-L)^{-1/2})\). The entire state and control trajectories remain bounded by Sokolov’s stabilization theorem.

Greater uncertainty seems to require more retained observations. Does the entire Xie--Guo positive regime admits a universal stabilizer with limited record budget? Such a budget, if it exists, could still be extremely large. Surprisingly, our result shows that six retained state-response records are enough. The magic number emerges from our scheme for deleting observations, which we specify below. Upon observing \(x_0\), set \(u_0 = 0\). At each time \(t\ge1\), insert the newly available pair \((x_{t-1},z_t)\) before deletion. If its location already occurs, replace the old pair at that location by the new pair.

\textbf{Deleting Scheme.} Each record contains an observed real state \(x_i\). We can therefore order the records by comparing these known values \(x_0,\ldots,x_{t-1}\). Remove a record when its two neighboring retained states are ``close enough'' compared with an adaptive resolution \(2 R_t\). \(R_t\) is updated recursively from the current observation and a few stored auxiliary scalars. If more than one records satisfy the deletion condition, delete the leftmost one. Repeat until no further deletion is possible. We call the collection of retained records a cache \(\mathcal C_t\). 

We retain the Xie--Guo stabilizer \eqref{eq:control}. At each step \(t\), we select the nearest state \(x_{i_t}\) from the cache \(\mathcal C_t\) and use its associated response. Break ties by choosing the leftmost one. The remaining task is to choose \(R_t\) so that storage stays uniformly bounded and stability is preserved. Luckily, \(R_t\) has a simple choice.

\textbf{Adaptive Resolution.} Define the observed state-range width by
\begin{equation}\label{eq:range_width}
B_t:=M_t-m_t
\end{equation}
and define its largest past one-step increase by
\begin{equation}\label{eq:largest_expansion}
H_t:=\max_{1\le s\le t}(B_s-B_{s-1}),
\qquad H_0:=0.
\end{equation}
Choose \(R_t\) by
\begin{equation}\label{eq:threshold}
R_t=\max\left\{H_t,\frac{B_t}{6}\right\}.
\end{equation}

We explain why at most six records remain. Suppose seven retained locations satisfy \(p_1<\cdots<p_7\). Since no further deletion is possible,
\begin{equation}
p_3-p_1>2R_t,\qquad
p_5-p_3>2R_t,\qquad
p_7-p_5>2R_t.
\end{equation}
Adding these inequalities gives \(p_7-p_1>6R_t\). Thus, \(B_t > 6R_t\). This contradicts the choice in \eqref{eq:threshold}.

Neighboring retained states are at most \(2R_t\) apart. Thus, any point between the extreme retained states is within distance \(R_t\) of a retained state. If \(x_t\) lies outside this interval, its distance to the nearest retained endpoint is \(B_t-B_{t-1}\le H_t\le R_t\). Consequently,
\begin{equation}\label{eq:resolution_guarantee}
d_t:=\min_{(p,z)\in\mathcal C_t}|x_t-p|\le R_t.
\end{equation}
The cache therefore provides an observation of the current state within distance \(R_t\). This justifies calling \(R_t\) the adaptive resolution.

\begin{figure*}[t]
    \centering
    \includegraphics[width=0.8\textwidth]{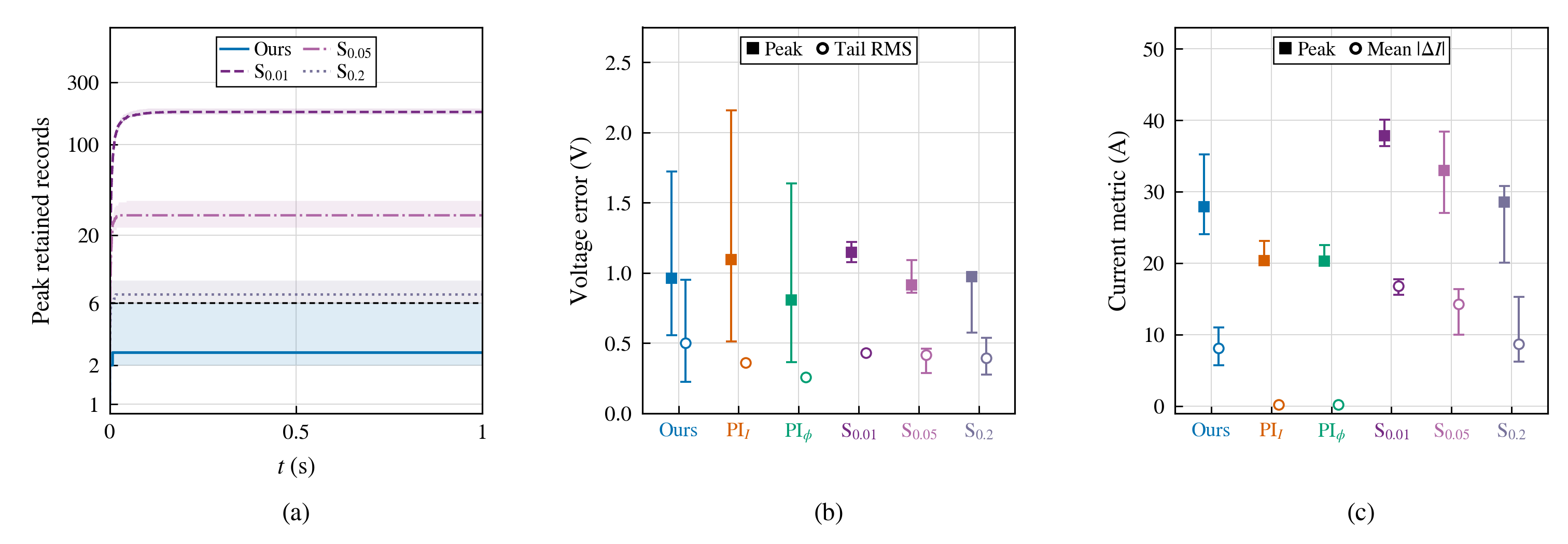}
    \caption{Comparison on the DAB DC-DC converter benchmark:
    (a) running maximum number of retained records (log scale);
    (b) peak and RMS voltage errors;
    (c) peak averaged-current command and mean absolute
    command change per sample.
    Curves and markers show medians over four initial phases;
    shading and whiskers show their ranges.
    $\mathrm{PI}_I$ and $\mathrm{PI}_\phi$ denote the
    linearization and feedback-only PI controllers, respectively.
    $S_\varepsilon$ denotes Sokolov's scheme with
    interval width $\varepsilon$ V.}
    \label{fig:comparison}
\end{figure*}

\textbf{Universal Stabilization with Six Records.} Now we prove the main result of this letter. Fix any \(f\) with \(\|f\|_{\mathrm{GL}}<L^\star\). Choose \(L\) such that \(\|f\|_{\mathrm{GL}}<L<L^\star\). By the definition of the generalized Lipschitz norm \eqref{eq:GL}, there exists a finite (unknown) constant \(c\ge0\) such that
\begin{equation}\label{eq:analysis_lipschitz}
|f(x)-f(y)|\le L|x-y|+c,
\qquad x,y\in\mathbb R.
\end{equation}
For a detailed discussion of this relation, see \cite[Remark~2.3 and Appendix~A]{XieGuo2000}. Set \(D:=c+2W\) for brevity. At each step \(t \geq 1\), the next state lies within distance \(LR_t + D\) of the midpoint of the observed state range \([m_t,M_t]\), because
\begin{equation}\label{eq:midpoint_bound}
\begin{aligned}
\left|x_{t+1}-\frac{m_t+M_t}{2}\right|
&=\left|f(x_t)+u_t+w_{t+1}-\frac{m_t+M_t}{2}\right|\\
&=\left|f(x_t)-z_{i_t}+w_{t+1}\right|\\
&=\left|f(x_t)-f(x_{i_t})+w_{t+1}-w_{i_t+1}\right|\\
&\le |f(x_t)-f(x_{i_t})|
+|w_{t+1}|+|w_{i_t+1}|\\
&\le L|x_t-x_{i_t}|+c+W+W\\
&=Ld_t+D \le LR_t+D.
\end{aligned}
\end{equation}
The successive steps use \eqref{eq:plant}, \eqref{eq:control}, \eqref{eq:visit_pair}, the triangle inequality, \eqref{eq:analysis_lipschitz}, and \eqref{eq:resolution_guarantee}, in that order. We translate this distance bound into a bound on the growth of the observed state range \(\Delta B_{t+1}:=B_{t+1}-B_t\). Consider where \(x_{t+1}\) falls relative to \([m_t,M_t]\). If \(x_{t+1}>M_t\), the width increases by \(x_{t+1}-M_t\). Thus
\begin{equation}
\begin{aligned}
\Delta B_{t+1}
&=x_{t+1}-M_t\\
&=\left(x_{t+1}-\frac{m_t+M_t}{2}\right)-\frac{B_t}{2}\\
&\le LR_t+D-\frac{B_t}{2}.
\end{aligned}
\end{equation}
The case \(x_{t+1}<m_t\) is handled similarly. If \(x_{t+1}\) lies inside the interval, the width does not increase. Combining these cases gives
\begin{equation}
\Delta B_{t+1}\le
\left[LR_t+D-\frac{B_t}{2}\right]^+,
\end{equation}
where \([v]^+:=\max\{v,0\}\). Substituting \eqref{eq:threshold} gives
\begin{equation}
\Delta B_{t+1}\le
\max\left\{ \left[LH_t-\frac{B_t}{2}+D\right]^+,\left[D-\frac{3-L}{6}B_t\right]^+ \right\}.
\end{equation}
Using this growth bound, we show that \(B_t\) is bounded. Suppose that \(B_t\) is unbounded. Since \(B_t\) never decreases, it must tend to infinity. Because \(L<L^\star<3\), the second term becomes zero once \(B_t\ge6D/(3-L)\). Therefore, for all sufficiently large \(t\),
\begin{equation}\label{eq:range_growth_criterion}
\Delta B_{t+1}
\le\left[LH_t-\frac{B_t}{2}+D\right]^+.
\end{equation}
Recall that \(B_t=\sum_{s=1}^t\Delta B_s\) and \(H_t=\max_{1\le s\le t}\Delta B_s\). By \cite[Lemma~3.3]{XieGuo2000}, \eqref{eq:range_growth_criterion} implies \(\sum_{t=1}^{\infty}\Delta B_t<\infty\), i.e. \(B_t\) is bounded, which leads to a contradiction. Consequently, \(B_\infty:=\sup_{t\ge0}B_t<\infty\). Since \(x_0,x_t\) and the midpoint lie in \([m_t,M_t]\), both \(x_t\) and the midpoint are within distance \(B_t\) of \(x_0\). The triangle inequality therefore gives
\begin{align}
|x_t|
&\le |x_0|+B_\infty,
\label{eq:state_bound}\\
\left|\frac{m_t+M_t}{2}\right|
&\le |x_0|+B_\infty.
\label{eq:midpoint_absolute_bound}
\end{align}
In particular, \eqref{eq:state_bound} shows that the state sequence \((x_t)_{t\ge0}\) is bounded. For every observation record \((x_i,z_i)\), we have
\begin{equation}\label{eq:response_bound}
\begin{aligned}
|z_i|
&=|f(x_i)+w_{i+1}|\le |f(x_i)|+|w_{i+1}|\\
&\le |f(x_0)|+|f(x_i)-f(x_0)|+W\\
&\le |f(x_0)|+L|x_i-x_0|+c+W\\
&\le |f(x_0)|+LB_\infty+c+W.
\end{aligned}
\end{equation}
Therefore, the entire input sequence is also bounded:
\begin{equation}
\begin{aligned}
|u_t|
&\le |z_{i_t}|+\left|\frac{m_t+M_t}{2}\right|\\
&\le |f(x_0)|+(L+1)B_\infty+c+W+|x_0|
<\infty.
\end{aligned}
\end{equation}

\textbf{Numerical Experiment.} Dual-active-bridge (DAB) DC--DC converters are widely used in energy storage and DC power systems, where voltage regulation under fluctuating loads is essential. We consider a DAB converter supplying a \(1.6\)-kW single-phase load with \(120\)-Hz power ripple \cite{QinKimball2014}. The input and reference voltages are \(400\) and \(160\) V, respectively. We use \(C=1\) mF, \(L=70\,\mu\mathrm H\), a transformer turns ratio of \(2{:}1\), and \(20\)-kHz sampling \cite{Shao2022}. In the Euler-discretized voltage dynamics, the controlled current produces an additive control input, while the load-current ripple enters as a bounded additive disturbance. The resulting plant lies in the Xie--Guo positive regime. We implement the Xie--Guo stabilizer using our six-record scheme and Sokolov's scheme with spatial resolutions of \(0.01\), \(0.05\), and \(0.2\) V. We also compare against the feedback-only and linearization-based PI controllers in \cite{Shao2022}, using their reported gains without retuning. Each \(1\)-s simulation uses fourth-order Runge--Kutta integration with a \(50\,\mu\mathrm s\) sampling period. Our scheme retains only \(2\)--\(6\) records, compared with \(171\)--\(190\) and \(23\)--\(37\) records for Sokolov's scheme at \(0.01\)-V and \(0.05\)-V resolutions, respectively. Moreover, its peak voltage deviation from the \(160\)-V reference is comparable to those of the two PI controllers, demonstrating effective voltage regulation with a small observation budget.

\textbf{Conclusion.} Classical adaptive control often relies on finite-dimensional recursive updates, but its guarantees typically concern structured uncertainty classes \cite{KarafyllisKrstic2018}. We investigate whether this principle survives at the sharp limit of universal stabilizability. In this letter, we show that just six observation records and a few auxiliary summaries suffice for stabilization throughout the entire Xie--Guo positive regime. Is six the smallest sufficient record budget? More fundamentally, how far does bounded retention extend beyond the geometry of generalized Lipschitz uncertainty? Which universally stabilizable uncertainty classes admit stabilization with bounded retention, and which inherently require growing memory? Answering these questions may reveal how the acquisition, organization, and use of knowledge shape the fundamental limits of control under uncertainty.

\textbf{Acknowledgments:} This work was supported by the National
Natural Science Foundation of China under Grant No. 12288201.


\end{document}